\documentclass[12pt]{amsart}

\usepackage{amssymb}

\newtheorem{theorem}{Theorem}[section]
\newtheorem{lemma}[theorem]{Lemma}
\newtheorem{corollary}[theorem]{Corollary}

\theoremstyle{definition}
\newtheorem{definition}[theorem]{Definition}
\newtheorem{example}[theorem]{Example}

\theoremstyle{remark}

\begin{document}

\title{Antichains of Monomial Ideals are Finite}

\author{Diane Maclagan}

\address{Department of Mathematics \\
University of California, Berkeley\\
Berkeley, CA 94720}

\email{maclagan@math.berkeley.edu}

\date{\today}

\begin{abstract}
The main result of this paper is that all antichains are finite in the
poset of monomial ideals in a polynomial ring, ordered by
inclusion. We present several corollaries of this result, both simpler
proofs of results already in the literature and new results.  One
natural generalization to more abstract posets is shown to be false.
\end{abstract}

\maketitle

\section{Introduction}

Throughout this paper, $S=k[x_1,\ldots ,x_n]$, where $k$ is
a field. Our main result is the following theorem:

\begin{theorem}\label{finiteantichains}
Let $\mathcal I$ be an infinite collection of monomial ideals in a
polynomial ring.  Then there are two ideals $I$, $J \in \mathcal I$
with $I \subseteq J$.
\end{theorem}

Although the statement may appear to be purely algebraic, monomial
ideals are highly combinatorial objects.  In particular, the above theorem
can be restated as follows:

\begin{theorem} \label{finiteposetform}
  Let $L$ be the poset of dual order ideals of the poset $\mathbb
N^n$, ordered by containment.  Then $L$ contains no infinite
antichains.
\end{theorem}

A special case of interest is Young's lattice, which consists of the
set of all partitions ordered by containment of Ferrers diagrams.
Noting that a partition can be considered to be a finite order ideal
in $\mathbb N^2$, we consider the generalized Young's lattice of finite
order ideals in $\mathbb N^n$ ordered by inclusion.

\begin{theorem} \label{young}
All antichains in the generalized Young's lattice are finite.
\end{theorem}

In the next section we give some corollaries of Theorem
\ref{finiteantichains}.  Some of the corollaries have appeared in the
literature before, but Theorem \ref{finiteantichains} allows us to
simplify the original proofs, and provides a common framework for
finiteness results involving monomial ideals.  In Section 3 we give an
application to SAGBI bases which was the motivating example for this
paper.  In Section 4 we outline an example  which
shows that one natural generalization to more abstract posets is false,
and lastly in Section 5 we give a proof of the theorem.

\section{Corollaries}

In this section we give several corollaries of Theorem \ref{finiteantichains}.

The first corollary is a new proof of a basic result in computational
algebra.  A fundamental notion in Gr\"obner basis theory is that of an
{\em initial ideal} of an ideal in a polynomial ring $S$.  Given a {\em
term order} $\prec$ (a total order on monomials in $S$ satisfying certain
conditions), we define the {\em initial term} of a polynomial to be
the largest monomial with respect to $\prec$ occurring in the
polynomial.  The initial ideal $in_{\prec}(I)$ of $I$ with respect to
$\prec$ is the monomial ideal generated by the initial terms of all
polynomials in $I$.  The following theorem appears in \cite{BM} and
\cite{MR}, and is well known.

\begin{corollary} \label{inideals}
For a given ideal $I \in S$ there are only finitely many distinct
 initial ideals $in_{\prec}(I)$.
\end{corollary}

\begin{proof}
The monomials of $S$ outside $in_{\prec}(I)$ form a $k$-basis for
$S/I$.  If there were infinitely many initial ideals then Theorem
\ref{finiteantichains} would give a proper inclusion of $k$-bases.
\end{proof}

Given an $\mathbb N^d$ grading on $S$, we can define the Hilbert series
of a homogeneous ideal by $$H_{S/I}(t)=\sum_{b \in \mathbb N^d} (\dim_k
(S/I)_b) t^b$$ where $t^b=\prod_{i=1}^d t_i^{b_i}$.

\begin{corollary}
There are finitely many monomial ideals with a given Hilbert series
with respect to  a given  grading. 
\end{corollary}

Theorem \ref{finiteposetform} is also true when $\mathbb N^n$ is
replaced by a finitely generated submonoid (such as the lattice points
inside a rational cone).

\begin{corollary}
Let $M$ be a finitely generated submonoid of $\mathbb N^n$.  Let
$R=k[M]=k[t^{a_1},\ldots,t^{a_d}]$ be its monoid algebra.  A monomial
ideal in $R$ is an ideal generated by elements of the form $t^b \in
R$ for some $b \in \mathbb N^n$.  Then in any infinite collection
$\mathcal I$ of monomial ideals in $R$ there are two, $I,J \in \mathcal I$,
such that $I \subseteq J$.
\end{corollary}
\begin{proof}
Consider the map $\phi:k[x_1,\ldots , x_d] \rightarrow R$ given by
$\phi:x_i \mapsto t^{a_i}$.  For a monomial ideal $I \subseteq R$, we
define $I_{\phi}= \langle x^a : \phi(x^a) \in I \rangle $.  Then $I_{\phi}
\subseteq J_{\phi} \Rightarrow I \subseteq J$, so the result follows
from applying Theorem \ref{finiteantichains} to the set $\mathcal I_\phi=\{I_{\phi}: I \in \mathcal I \}$.
\end{proof}

A similar corollary relates to $A$-graded algebras, where $A$ is a $d
\times n$ matrix with entries in $\mathbb N$.  An $A$-graded algebra
is a $k$-algebra $R$ generated by $x_1,x_2,\ldots ,x_n$ with an
$\mathbb N^d$ grading (given by $\deg x_i=a_i$, where $a_i$ is the
$i$th column of $A$) such that $\dim_k R_b =1$ whenever $b \in \mathbb
N A$ (the image of the map $\pi:\mathbb N^n \rightarrow \mathbb N^d$
given by $\pi:\omega \mapsto A \omega$) and equals $0$ otherwise.  See
\cite[Chapter 10]{GBCP} for details of $A$-graded algebras.

\begin{corollary} \label{agas}
Let $R$ be an $A$-graded algebra.  Let $\mathcal I$ be an infinite
collection of ideals of $R$ which are homogeneous with respect to the
$A$-grading.  Then there are two ideals, $I,J \in \mathcal I$ such that $I
\subseteq J$.
\end{corollary}
\begin{proof}
$R$ is isomorphic to $S/I$ for some binomial ideal $I$.  Any element
of $S/I$ which is homogeneous with respect to the $\mathbb N^d$
grading can be written as $m+I$ where $m$ is some monomial in $S$, so
homogeneous ideals of $R$ lift to monomial ideals in $S$.  Containment
in $S$ implies containment in $R$, so the result follows.
\end{proof}

A trivial example of an $A$-graded algebra is $k[x_1,\ldots ,x_n]$ with
$A$ the $n \times n$ identity matrix.  Then Corollary \ref{agas}
reduces to Theorem \ref{finiteantichains}.

\section{Application to SAGBI bases}

Let $T=R[c_1x^{a_1},\ldots ,c_nx^{a_n}]$ be a monomial subalgebra of
$R[x_1,\ldots ,x_d]$, where $R$ is a Principal Ideal Domain.  A strong
SAGBI (Subalgebra Analogue to Gr\"obner Bases for Ideals) basis for
$T$ is a collection $\{k_1x^{b_1},\ldots k_mx^{b_m} \}$ such that any
element $cx^l \in T$ can be written as $cx^l=r\prod_{i=1}^m
(k_ix^{b_i})^{\phi_i}$ for some $\phi \in \mathbb N^m$ and $r \in R$.

\begin{definition} \label{oldatomic}
Given a matrix $A \in \mathbb N^{d \times n}$, we define a map
$\pi:\mathbb N^n \rightarrow \mathbb N^d$ by $\pi: y \mapsto A
y$.  Let $\mathbb N A \subseteq \mathbb N^d$ be the image of
$\pi$.  For $b \in \mathbb N A$ let $P_b=conv(\pi^{-1}(b))$.  Since
$\pi^{-1}(b)$ is a finite set, this is a convex polytope.  We call
$P_b$ the fiber of $A$ over $b$.  A fiber over $b$ is {\em atomic} if
there do not exist $b_1,b_2 \in \mathbb N^d$ with $b_1+b_2=b$ such
that $P_b=P_{b_1}+P_{b_2}$, where the addition is Minkowski sum.
\end{definition}

Atomic fibers were defined by Adams et al.\  in \cite{AHLM}, where they
proved that there are only a finite number of atomic fibers for a
given matrix $A$.  They used this result to construct a finite strong
SAGBI basis as follows:

\begin{theorem}[Adams et. al.\ \cite{AHLM}]
Let $T=R[c_1x^{a_1},\ldots ,c_nx^{a_n}]$.  Let $A=[a_1,\ldots ,a_n]$ be
the $d \times n$ matrix with columns the $a_i$.  Then $\{k_{b}x^b:P_b$
is an atomic fiber of $A \}$ is a strong SAGBI basis for $T$, where
$k_b=\gcd(\{c^u=c_1^{u_1}\ldots c_l^{u_l} : u=(u_1,\ldots,u_n) \in
\pi^{-1}(b)\})$.
\end{theorem}

  The proof of the finiteness result in \cite{AHLM} was constructive but
complicated, using convex geometry techniques.  Theorem
\ref{finiteantichains} gives a much simpler, though non constructive,
proof of this result.

\begin{corollary} \label{atomic}
For a given matrix $A \in \mathbb N^{d \times n}$, there are only a finite number of atomic fibers.
\end{corollary}
\begin{proof}
For $b \in \mathbb N A$, let $I_b= \langle x^u : Au=b$ and $u$ is a
vertex of $P_b \rangle$.  Then the fiber over $b$ is atomic if and
only if $I_b$ is not contained in any $I_{b^{\prime}}$ for $b \neq
b^{\prime}$.  If there were an infinite number of atomic fibers, then $\{I_b : P_b \text{ atomic}\}$ would be an infinite antichain of monomial ideals, contradicting Theorem \ref{finiteantichains}.
\end{proof}

Corollary \ref{atomic} can be generalized as follows:

\begin{definition}
Let $M$ be a monomial ideal of $S$, and $A \in \mathbb N^{d \times n}$
a matrix.  Then the $(M,A)$ fiber over $b \in \mathbb N A$ is the set
$\{ u : Au=b$ and $x^u \notin M \}$.  A $(M,A)$ fiber over $b$ is
atomic if there do not exist $b_1,b_2 \in \mathbb N A$ with $b_1+b_2=b$
such that for all $u$ in the $(M,A)$ fiber over $b$ there are $u_1,
u_2$ in the $(M,A)$ fibers over $b_1, b_2$ respectively such that
$u=u_1+u_2$.
\end{definition}

To see that this definition is a generalization of an earlier one, we first need another definition.

\begin{definition} \label{newatomic}
Given a matrix $A \in \mathbb N^{d \times n}$, we define its vertex
ideal, $V_A$ by $$V_A= \bigcap_{\prec} in_{\prec}(I_A)$$ where the intersection is over all term orders $\prec$, and $I_A$ is the toric ideal corresponding to $A$ (see \cite{GBCP} for details on toric ideals).  
\end{definition}

Note that this is a finite intersection by Corollary \ref{inideals}.
Since the standard monomial of $A$-degree $b$ of an initial ideal of a
toric ideal corresponds to a vertex of $P_b$, and each vertex of 
$P_b$ is standard for some initial ideal, the set of standard
monomials of $V_A$ is exactly $\{x^u : u$ is a vertex of $P_{Au} \}$.
Thus Definition \ref{oldatomic} is Definition \ref{newatomic} with
$M=V_A$.

\begin{corollary} \label{atomictwo}
There are only finitely many atomic $(M,A)$ fibers for given $M$ and $A$.
\end{corollary}

The proof is the same as for Corollary \ref{atomic}.  Of particular
interest is the case $M=(0)$.  In that case,  being atomic
corresponds to the nonexistence of $b_1,b_2$ such that each lattice point
in $\pi^{-1}(b)$ is a sum of lattice points in $\pi^{-1}(b_1)$ and
$\pi^{-1}(b_2)$, as opposed to the original definition, where only the
vertices need be sums of lattice points in the two smaller fibers.
This is a strictly stronger requirement.  The following example shows
that a fiber can be atomic with respect to this stronger definition
without being atomic in the original sense.

\begin{example}
Let $A$ be the following matrix:
$$\left( \begin{array}{cccccc}
1 & 1& 1& 0 & 0 & 0\\
0 & 3 & 2 & 1 & 0 & 0\\
5 & 0 & 2 & 0 & 1 & 0\\
0 & 2 & 1 & 0 & 0 & 1\\
\end{array} \right)$$

Let $b_1=(1,3,5,2)^T$, and $b_2=(5,10,10,6)^T$.  We have
\begin{eqnarray*}
\pi^{-1}(b_1) & = \{ &(1, 0, 0, 3,0,2)^T,\\
& & (0,1,0,0,5,0)^T,\\
& &(0,0,1,1,3,1)^T \}\\
\end{eqnarray*}
 and
\begin{eqnarray*}
\pi^{-1}(b_2) & =\{ &(0,0,5,0,0,1)^T,\\
& &(1,2,2,0,1,0)^T,\\ 
& & (2,3,0,1,0,0)^T \}\\
\end{eqnarray*}
 Now $P_{b_1+b_2}=P_{b_1}+P_{b_2}$, so $b_1+b_2$ is not atomic in the
first sense.  However $(1,1,4,2,2,2)^T \in \pi^{-1}(b_1+b_2)$, but
cannot be written as the sum of lattice points in $P_{b_1}$ and
$P_{b_2}$. This example is based on an example of Oda \cite{oda} for
lattice polytopes.
 
\end{example}

\section{Poset Formulation}

{}From the second formulation of the theorem, it is natural to suspect
that this is in fact a general theorem about posets.  Two properties
of the poset $\mathbb N^n$ which lend themselves to finiteness results
are that $\mathbb N^n$ has no infinite antichains, and satisfies the
descending chain condition.  The following example consists of a poset
which has no infinite antichains or infinite descending chains  such
that the poset of dual order ideals under containment contains an
infinite antichain.  This example appears in \cite{DPR}, but was
discovered independently by George Bergman, from whom I learned it.

\begin{example} {\em (\cite{DPR}, G. Bergman)} Let $X$ be the set $\{ (i,j) : i,j
\in \mathbb N, i<j \}$.  Set $(i,j) \prec (i^{\prime},j^{\prime})$ if
and only if $j<j^{\prime}$ and either $i=i^{\prime}$ or
$j<i^{\prime}$.  It is straightforward to check that $X$ is a
partially ordered set.

Note that any chain descending from $(i,j)$ can have at most $j-1$
members less than $(i,j)$, so there are no infinite descending chains
of elements of $X$.  To see that all antichains in $X$ are finite,
suppose $Y$ is an antichain in $X$, and let $j_0$ be the smallest $j$
such that $(i,j) \in Y$, occurring in the pair $(i_0,j_0) \in Y$.
Then $(i,j) \in Y$ implies $i \leq j_0$, as otherwise $j_0 < j$, and
then  $(i_0,j_0) \prec (i,j)$.  If there are two pairs $(i,j),
(i,j^{\prime}) \in Y$, with $j < j^{\prime}$ then $(i,j) \prec
(i,j^{\prime})$, so there is only one pair of the form $(i,j) \in Y$
for each value of $i$.  But this means there are at most $j_0+1$
elements in $Y$, so all antichains in $X$ are finite.

Because there are no infinite descending chains or infinite antichains
each dual order ideal in $X$ can be represented by its finite
antichain of minimal elements.  One dual order ideal is contained in
another exactly when each element of the finite antichain of minimal
elements of the first dual order ideal is greater than some element of
the finite antichain of minimal elements of the second.

For fixed $l >0$, let $S_l=\{ (k,l) : k < l \} \subseteq X$.  Then $S_l$ is
the finite antichain of minimal elements of a dual order ideal of $X$.
Suppose the dual order ideal determined by $S_{l_2}$ is contained in the
one determined by $S_{l_1}$.  From above, we must have $l_1 < l_2$.  But then
there is no element of $S_{l_1}$ less than $(l_1,l_2) \in S_{l_2}$, a
contradiction.  So the $S_l$ form an infinite antichain of dual order
ideals of $X$.
\end{example}

Theorem \ref{finiteposetform} can, however, be generalized in the following way:

\begin{theorem} (Farley, Schmidt) \cite{Farley}
Let $P$ and $Q$ be two posets with no infinite antichains that satisfy
the descending chain condition.  If the posets of dual order ideals of
$P$ and of $Q$, ordered by inclusion, have no infinite antichains, then
the same is true for the poset of dual order ideals of $P \times Q$.
\end{theorem}

\section{Proof of the Main Theorem}

In this section we outline the  proof of Theorem \ref{finiteantichains}.

We first prove Theorem \ref{young}.  The generalized Young's lattice
is isomorphic to the poset of artinian monomial ideals under
inclusion, via the map taking an order ideal to its complement, so we
prove the theorem in that setting.

\begin{lemma}\label{young2}
Let $\mathcal I$ be an infinite collection of artinian monomial ideals
(primary to the maximal ideal).  Then there are two ideals, $I,J \in
\mathcal I$ such that $I \subseteq J$.
\end{lemma}

\begin{proof}
Suppose $\mathcal I$ consists of an infinite number of artinian
monomial ideals, which are noncomparable with respect to inclusion.
Choose $I_1 \in \mathcal I$.  Since $I \not \subseteq I_1$ for $I \in
\mathcal I \setminus \{I_1\}$, each $I \in \mathcal I \setminus
\{I_1\}$ contains some of the finite number of standard monomials of
$I_1$.  There are thus an infinite number of ideals in $\mathcal I$
which contain the same set of standard monomials of $I_1$.  Call this
collection $\mathcal I_1$.  Let $J_1$ be the intersection of the
ideals in $\mathcal I_1$.  We will now build a strictly ascending
chain of monomial ideals.  Suppose $\mathcal I_k$ and $J_k$ have been
chosen.  Choose an ideal $I_{k+1} \in \mathcal I_k$.  We can again
find an infinite collection of ideals in $\mathcal I_k$ which have the
same non-trivial intersection with the standard monomials of
$I_{k+1}$.  Let $\mathcal I_{k+1}$ be this collection, and let
$J_{k+1}$ be the intersection of the ideals in $\mathcal I_{k+1}$.  We
have $J_{k+1} \supsetneq J_k$, since $J_{k+1}$ contains some standard
monomials of $I_{k+1}$, so in this fashion we get an infinite
ascending chain of monomial ideals in $S$, which is impossible.
\end{proof}

\begin{corollary} \label{primarychain}
Let $\mathcal I$ be an infinite collection of artinian monomial
ideals. Then there is an infinite chain $I_1 \supsetneq I_2 \supsetneq
\ldots$ of ideals in $\mathcal I$.
\end{corollary}
\begin{proof}
Since $S$ is Noetherian, $\mathcal I$ contains maximal ideals.  There
are only finitely many maximal ideals by Lemma \ref{young2}, so set
$I_1$ to be a maximal ideal in $\mathcal I$ containing an infinite
number of ideals of $\mathcal I$, and repeat, setting $\mathcal I= \{
I \in \mathcal I : I \subsetneq I_1\}$.
\end{proof}

\begin{proof}[Proof of Theorem \ref{finiteantichains}.]
Every associated prime of a monomial ideal is a monomial prime, of
which there are only a finite number.  We can thus restrict to an
infinite collection of $\mathcal I$ all of which have the same set of
associated primes, which we will also call $\mathcal I$.  Now for each
ideal in this set we find an irredundant primary decomposition,
writing the ideal as the intersection of monomial ideals primary to an
associated prime in such a way that each associated prime is used only
once.  Let $I_{\tau}$ be the primary component of $I$ primary to the
monomial prime $P_{\tau}=\langle x_i : i \not \in \tau \rangle$, where
$\tau \subseteq [n]$.  For a fixed $\tau$ either $\{ I_{\tau}: I \in
\mathcal I \}$ is finite, so there is an infinite number of $I \in
\mathcal I$ with the same $I_{\tau}$, or we can apply Corollary
\ref{primarychain} to the polynomial ring $k[x_i : i \not \in \tau]$.
In either case we get an infinite collection $\mathcal I_{\tau}=\{ I_k
: k \geq 1 \}$ of ideals in $\mathcal I$ such that ${I_1}_{\tau}
\supseteq {I_2}_{\tau} \supseteq \ldots $, where the inclusions need
not be proper.  Since there are only a finite number of associated
primes, by appropriate restrictions we can find a sequence $\{ I_k : k
\geq 1 \}$ such that ${I_1}_{\tau} \supseteq {I_2}_{\tau} \supseteq
\ldots$ for each $\tau$ such that $P_{\tau}$ is an associated prime.
But since $I_k$ is the intersection of the ${I_k}_{\tau}$, where
$\tau$ ranges over all over associated primes $P_{\tau}$ of $I_k$,
this means that $I_1 \supsetneq I_2 \supsetneq \ldots$, where the
inclusions are proper, since the $I_k$ are all distinct.
\end{proof}

\section{Acknowledgments}
The proof of Theorem \ref{finiteantichains} was much improved from a previous version by discussion with Dave Bayer.

\end{document}